\documentclass[12pt,a4paper]{article}
\RequirePackage{etex}
\usepackage{amsfonts, amsmath, amssymb,amsthm,graphics,mathrsfs,verbatim, CJK}
\usepackage{makeidx}
\usepackage{palatino}
\usepackage{titletoc}
\usepackage{tocloft}
\usepackage{titlesec}
\usepackage{stmaryrd,color}
\usepackage[colorlinks=true]{hyperref}
\usepackage{graphicx}
\usepackage{pgfplots}
\usepackage{tikz}
\usepackage{pstricks,pst-node,multido,ifthen,calc}

\usetikzlibrary{arrows}
\usetikzlibrary{calc}

\pgfplotsset{compat=1.16}

\titleformat{\section}{\Large}{}{0.5em}{}
\makeindex

\normalsize

\newtheorem{mthm}{Theorem}[section]
\newtheorem{mylem}[mthm]{Lemma}
\newtheorem{myprn}[mthm]{Proposition}
\newtheorem{mycor}[mthm]{Corollary}
\newtheorem{mydef}[mthm]{Definition}
\newtheorem{myrem}[mthm]{Remark}
\newtheorem{mycon}[mthm]{Construction}
\newtheorem{myeg} [mthm]{Example}
\newtheorem{myque} [mthm]{Question}
\newtheorem{myalg} [mthm]{Algorithm}
\newtheorem{myconj} [mthm]{Conjecture}
\newenvironment{thm}{\begin{mthm}}{\end{mthm}}

\def \Lemma #1 {\vs{2mm}\nin {\bf Lemma #1} }
\def \Prop #1 {\vs{2mm}\nin {\bf Proposition #1} }
\def \Th #1 {\vs{2mm}\nin {\bf Theorem #1} }
\def \Cor #1 {\vs{2mm}\nin {\bf Corollary #1} }
\def \Proof{\vs{2mm}\nin {\bf Proof.~}}
\def \Def #1 {\vs{2mm}\nin {\bf Definition #1} }
\def \part #1 {\hfil\break\hglue 12pt {\rm (#1)~}}

\def\a{\alpha}

\def\D{\Delta}
\def\hs{\hspace}

\def\eset{\emptyset}

\def\fs{\footnotesize}

\def\G{\Gamma}

\def\nin{\noindent}

\def \qed {~\vrule height6pt width 6pt depth 0pt}

\def\LRA{\longrightarrow}

\def\si{\sigma}

\def\seq{\subseteq}

\def\vs{\vspace*}

\def\CM{ Cohen-Macaulay }
\def\iff{ if and only if }

\def\vd{vertex-decomposable}

\makeatletter
\newcommand{\eqnum}{\refstepcounter{equation}\textup{\tagform@{\theequation}}}
\makeatother

\begin{document}
\title{
\bf\LARGE  Simplicial complexes which are minimal Cohen-Macaulay\thanks{This research was supported by the Natural Science
Foundation of Shanghai (No. 19ZR1424100), and the National Natural Science Foundation of China (No. 11971338) } }

\author{ Yanyan Wang\thanks{yanyanwang0411@sjtu.edu.cn},\,\,\, Tongsuo Wu\thanks{Corresponding author. {\small
tswu@sjtu.edu.cn}}\\
{\fs School of Mathematical Sciences, Shanghai Jiao Tong University}\\}

\date{}
\maketitle

\begin{center}
\begin{minipage}{12cm}

\vs{3mm}\nin{\small\bf Abstract} {\fs Let $\D$ be a $(d-1)$-dimensional pure  $f$-simplicial complex  over vertex set $[n]$. In this paper, it is proved that $\D$ being minimal CM implies $d\ge 3$ and  $n=2d$. It is also indicated that the recent work of \cite{Dao2020} implies that shellable condition on a pure simplicial complex $\D$ is identical with  existence of a full series of CM subcomplexes of $\D$.}

\vs{3mm}\nin {\small\bf Key Words} {\small Simplicial complex; Cohen-Macaulay; shellable;  shelled over; $f$-simplicial complex }

\vs{3mm}\nin {\small 2020 AMS Classification:} {\small Primary: 13F55; 05E45; Secondary: 13H10; 55U10.}

\end{minipage}
\end{center}

\section{1. Introduction}

For a natural number $n$ and $1<d< n$, let  $$[n]=:\{1,2,\ldots,n\},\,\,[n]_d=:\{A\in 2^{[n]}\mid\,\,\, |A|=d\},$$
where $2^{[n]}$ is the power set of $[n]$. A simplicial complex $\D$ over a vertex set $[n]$ is a subset of $2^{[n]}$, which has the hereditary property under inclusion and is such that $\{i\}\in \D$ holds for all $i\in [n]$.  Recall that a facet of $\D$ is a maximal element with respect to inclusion, and the facet set of $\D$ is denoted as $\mathcal F(\D)$. The dimension ${\rm dim}\,\D$ of $\D$ is the maximal number $|F|-1$, where $F\in \mathcal F(\D)$ runs over all facets of $\D$. If ${\rm dim}\,\D$ equals to $|F|-1$ for each facet $F$, then $\D$ is said to be {\em pure}. Let $\D^{(i)}=\{F\in \D\mid |F|\leq i+1\}$ be the $i$'th skeleton of $\D$.

Cohen-Macaulay (abbreviated as CM) property is one of the central research topics in commutative algebra and the rich and deep homological achievements have fruitful applications in combinatorial aspects of commutative rings (\cite{Stanley,Eisenbud,Bruns and Herzog,Villarreal,MS,Herzog and Hibi}). In combinatorial commutative algebra, shellable and pure simplicial complexes are the main source of CM simplicial complexes. In a most recent work \cite{Dao2020}, Dao, Doolittle and  Lyle discovered a new important combinatorial property of a CM simplicial complex $\D$, i.e., {\em $\D_F\cap \langle F\rangle $ is pure of dimension $|F|-2$} for any facet $F$ of $\D$, where $\mathcal F(\D_F)=\mathcal F(\D)\smallsetminus \{F\}$. Based on the property, the notion of a minimal CM simplicial complex $\D$ is introduced and studied. To be more precisely, $\D$ is called minimal CM if $\D$ is CM but no $\D_F$ is CM for any facet $F$ of $\D$. Acyclic behavior of a minimal CM $\D$ is studied and, sufficient conditions are provided for a complex to be minimal CM . Many interesting examples of minimal CM complexes are also exhibited.   Recall also from Zheng \cite{ZhengXX} the other important combinatorial property of a CM simplicial complex, i.e., {\em CM simplicial complexes  are connected in codimension one}, i.e., for any distinct facets $F$ and $G$,  there is a sequence $F=F_0, F_1, \ldots, F_r=G$ of facets such that $|F_i\cap F_{i+1}|=|F_{i+1}|-1$ for all $i=0, 1, \ldots, r-1$.

In this paper, we use  \cite[Lemma 3.1]{Dao2020} to study the exact relation of  shellable and CM properties for a pure complex $\D$, and we study the condition for a  minimal CM $f$-simplicial complex to be acyclic. In Section $2$, we recall some work of \cite{Dao2020} and, give a brief survey on $f$-simplicial complexes. In Section $3$, we first indicate that shellable condition on a pure complex $\D$ is identical with  CM properties of a full series of subcomplexes of $\D$, and then we use this observation to construct nontrivial examples of pure shellable complexes by taking advantage of CoCoA in an algorithmic approach, after applying Eagon-Reiner theorem. In Section $4$, we compute the dimension of the subspace ${\rm ker}(\partial_r)$ in a reduced chain complex of the simplex $\langle [n]\rangle$, and apply it to deduce that a minimal CM $f$-simplicial complex exists in $[n]_d$ implies $n=2d$.

\section{2. Preliminaries}

 For a $(d-1)$-dimensional simplicial complex $\D$, there is a related chain complex of $\mathfrak{K}$-spaces:
$$\mathcal C:\,\,\, 0\LRA C_{d-1}\overset{\partial_{d-1}}{\LRA} C_{d-2}\overset{\partial_{d-2}}{\LRA}\cdots \LRA C_{1}\overset{\partial_1}{\LRA} C_0\LRA 0,$$
where $C_i$ is a free $\mathfrak{K}$-module with basis set $\{\si\in \D\mid |\si|=i+1\}$, while for any $1\le k_1<k_2<\cdots <k_{r+1}\leq n$,
$$\partial_r(k_1k_2\ldots k_{r+1})=\sum_{i=1}^{r+1}(-1)^{i-1}k_1\ldots k_{i-1}\hat{k_{i}}k_{i+1}\ldots k_{r+1}.$$ For each $i$, recall that ${\rm im}\,\partial_{i+1}\seq \ker\, \partial_{i} $ holds, and the quotient $\mathfrak{K}$-space $$\tilde{H}_i(\D)=:\ker\,\partial_i/{\rm im} \,\partial_{i+1}$$ is called the $i^{th}$ homology group of $\D$. If $\tilde{H}_i(\D)=0$ holds for all $i$, then $\D$ is said to be {\em acyclic}. Clearly, $\D$ is acyclic \iff the corresponding chain complex $\mathcal C$ is an exact sequence. Recall that a cone is always acyclic  (see, e.g., \cite{Villarreal}), where $\D$ is called a cone if there exists a vertex such that all facets contain it as an element, and note that $\dim\,\tilde{H}_0(\D)+1$ is the number of connected components of $\D$ (\cite[Proposition 5.2.3]{Villarreal}).

Throughout, let $\mathfrak K$ be a field and let $S=\mathfrak K[x_1,\ldots,x_n]$ be the polynomial ring over $\mathfrak K$. Throughout, unless otherwise specifically stated, let $\D$ be a $(d-1)$-dimensional pure simplicial complex with vertex set $[n]$, where $\D\not=\eset,\{\eset\}$ and $\D$ is not a simplex.

We first recall some work of \cite{Dao2020} on minimal Cohen-Macaulay simplicial complexes. For a facet $F$ of a simplicial complex $\D$, let
$$\D_F=\langle G\mid G\in \mathcal F(\D)\smallsetminus \{F\}\rangle.$$ $\D$ is called a {\em shelling move} of $\D_F$ if $\D_F\cap \langle F\rangle$ is pure of codimension $1$, i.e., $\D_F\cap \langle F\rangle$ is generated by some nonempty subset of $\partial F$.

\begin{mylem} ({\rm \cite[Lemma 3.1]{Dao2020}})\label{shelling move} If a simplicial complex $\D$  is CM with  $|\mathcal F(\D)|\ge 2$, then $\D$ is a shelling move of $\D_F$ for any facet $F$ of $\D$.
\end{mylem}

A CM simplicial complex $\D$ is called {\em minimal} CM, if either it is a simplex, or else $|\mathcal{F}(\D)|\ge 3$ and,  no $\D_F$ is CM for any facet $F$ of $\D$. Here the definition of minimal CM is slightly different from that of \cite{Dao2020}, since we are mainly interested in the nonempty simplicial complexes. Then, the shelling move property implies the following:

\begin{mthm} {\rm (\cite[Theorem 3.2]{Dao2020})}\label{shelled over main}  Let $\D$  be a $(d-1)$-dimensional CM simplicial complex, which is not minimal. Then there exists a minimal CM subcomplex $\G$ and a series of facets $F_j,\ldots, F_1$ of $\D$, such that each $\G\cup \langle F_i,\ldots, F_1\rangle$ is CM and, each $\G\cup \langle F_1,\ldots, F_{i+1}\rangle$ is a shelling move of $\G\cup \langle F_1,\ldots, F_i\rangle$.
\end{mthm}

In \cite{Dao2020}, $\D$ is said to be {\it shelled over} $\G$. Clearly, shelled over is a kind of generalization of shellable for a simplicial complex.

Next, we record the following result, which is needed in this paper:

\begin{mthm} {\rm (\cite[Theorem 3.4]{Dao2020})}\label{facetdepth} Let $\D$ be a simplicial complex  with dimension $d-1$, and assume $\tilde{H}_{d-1}(\Delta) \ne 0$. Then there is a maximal facet $F$ of $\Delta$ such that the following hold:

$(1)$
$\dim \tilde{H}_{i-1}(\Delta_{F})=\begin{cases} \dim \tilde{H}_{i-1}(\Delta), \,\,\, \mbox{if } 0 \le i<d \\ \dim \tilde{H}_{i-1}(\Delta)-1, \,\, \mbox{if } i=d \end{cases}. $

$(2)$
$f_{k-1}(\Delta_{F})=\begin{cases} f_{k-1}(\Delta), \,\,\, \mbox{if } 0 \le i<d \\ f_{k-1}(\Delta)-1, \,\,\, \mbox{if } i=d \end{cases}.$

$(3)$
${\rm depth}\, \Delta={\rm depth} \,\Delta_{F}.$

\end{mthm}

Surely, this theorem together with Reisner theorem imply that a minimal CM simplicial complex is acyclic.


\begin{mycor}\label{acyclic2} Let $\D$ be a minimal CM complex over vertex set $[n]$ with dimension $d-1$. Then
$\tilde{H}_{d-1}(\Delta_F)=0$ holds true for any facet $F$ of $\D$.
\end{mycor}

\index{simplicial complex!minimal CM}

\Proof This is an immediate consequence of Lemma \ref{shelling move}, Theorem \ref{facetdepth} and \cite[Theorem 25.1,  page 142]{Munkres}.
\quad\quad\qed

Now we give a brief survey on some related established results on $f$-simplicial complexes. For any square-free monomial ideal $I$ of $S$, let $G(I)$ be the set of minimal monomial generators, and let ${\rm sm}(I)$ be the set of square-free monomials. For the ideal $I$, recall that there exist two related simplicial complexes, i.e., the nonface simplicial complex
$\delta_{\mathcal{N}}(I)=:\{\,  F\in 2^{[\, n\, ]}\mid X_F\in {\rm sm}(S)\smallsetminus {\rm sm}(I) \,  \}$
of $I$ and the facet simplicial complex
$\delta_\mathcal{ F}(I)=:\langle\,  F\in 2^{[\, n\, ]}\mid X_F\in {\rm G}(I)\, \rangle$
of the clutter $G(I)$. If they possess  a same $f$-vector, then the  ideal $I$ is called an {\it $f$-ideal}. For a simplicial complex  $\D$,  if its facet ideal
$I(\D)=:\langle\,\{X_F\mid F\in \mathcal F(\D)\}\,\rangle$
is an $f$-ideal,  then $\D$ is called an {\it $f$-simplicial complex}. A graph $G$ is said to be an $f$-graph, if the edge ideal $I(G)$ is an $f$-ideal. Note that in defining an $f$-graph $G$, $G$ is regarded as a simplicial complex of dimension no more than $1$, although we do have $I(G)=I_{{\rm Ind}(G)}$, where ${\rm Ind}(G)$ is the independence simplicial complex of the graph $G$. Refer to \cite{2012Anwar,2014Anwar,GuoWuFideals2,2016Anwar,LiuGuoWu5} for further related studies.

For a simplicial complex $\D$ on the vertex set $[n]$, let
$$\D^c=:\langle\,\{ F\mid F^c\in \mathcal F (\D)\}\,\rangle\overset{i.e.}{=}\langle\, \{[n]\smallsetminus G\mid G\in \mathcal F (\D)\}\,\rangle.$$
The definition of an $f$-simplicial complex seems to be reasonable with hindsight, due to the following two theorems on $f$-ideals.

\begin{thm}\label{homogeneous $f$-Ideal} {\rm (\cite[Theorem 2.3]{GuoWuLiu})} Let $S=K[x_1,  \ldots,  x_n]$,  and let $I$ be a  square-free monomial ideal of $S$ with the minimal generating set $G(I)$, where all monomials of $ G(I)$ have a same homogeneous degree $d$. Then $I$ is an $f$-ideal if and only if,  the set $ G(I)$ is an LU-set and,  $| G(I)| = \frac{1}{2}\binom{n}{d}$ holds true.
\end{thm}

Note that $ G(I)$ is said to be an L-set (U-set, respectively) if the set of all degree $d-1$ factors of elements of $G(I)$ has exactly $\binom{n}{d-1}$ elements (respectively, the set of degree $d+1$ square-free monomials extended from elements of $G(I)$ has cardinality $\binom{n}{d+1}$). If $G(I)$ is both a U-set and an L-set, then $G(I)$ is an LU-set. With the bijection from $X^\a$ to $\{i\in [n]\mid i\in \a\}$ (e.g., $x_1x_3x_4\mapsto \{1,3,4\}$), one obtains the notion of an L-set (U-set, LU-set respectively) for the facet set $\mathcal F(\D)$.

Recall also the following recently discovered result:

\begin{thm} {\rm(\cite[Theorem 4.1]{BTuyl2018})}\label{Newton} $\D $ is an $f$-simplicial complex, if and only if $\D^c$ is an $f$-simplicial complex.
\end{thm}

Equivalently, a square-free monomial ideal $I$ of $S$ is an $f$-ideal if and only if the Newton complement dual ideal
$\hat{I}=\langle x_1x_2\cdots x_n/u\mid u\in G(I)\rangle$
of $I$ is an $f$-ideal.

It is clear that  Theorem \ref{Newton} follows easily from Theorem \ref{homogeneous $f$-Ideal} for a {\it pure} simplicial complex $\D$.

Recall that for an  $f$-graph $G$, it is proved that the complement graph $\overline{G}$ is bipartite, thus ${\rm Ind(G)}=\overline{G}$. Recall that all $f$-graphs  are pure shellable as a graph, i.e., the independence complex ${\rm Ind(G)}$ is pure and shellable (\cite[Theorem 6.5]{GuoWuLiu}), while the definition of an $f$-graph is actually an $f$-simplicial complex of dimension less than or equal to $1$. Thus \cite[Theorem 6.5]{GuoWuLiu} may be re-stated as the following:

{\em If $\D$ is an $f$-simplicial complex of dimension less than or equal to $1$, then the homogeneous complement simplicial complex $\D'=:\langle\, [n]_2\smallsetminus \mathcal F(\D)\,\rangle$ is pure shellable.}

Based on this observation, it is natural to ask the following question:

\begin{myque}\label{Homogeneous complement question} For a pure $f$-simplicial complex $\D$ of dimension $d-1$, is the homogeneous complement simplicial complex $\D'=:\langle\,\si\mid\si\in [n]_d\smallsetminus \mathcal F(\D)\,\rangle$ of $\D$  shellable?

\end{myque}

We do not know counterexample in $[5]_3$ and in $[6]_3$. But it fails in $[8]_4$. We will give a negative answer in Example \ref{[8]4}.

We remark that there exist a lot of pure $f$-simplicial complexes which are not CM when $d-1\ge 2$.

Finally, we claim that there exist $f$-simplicial complexes which are minimal CM:

\begin{myeg} {\rm (\cite{BW1996, LiuGuoWu5})}\label{f-idealNonShellable} {\rm
 Consider a simplicial complex $\D$ with facet set
$$\mathcal F(\D)=\{123,  125,  136,  145,  146,  234,  246,  256,  345,  356 \}$$
constructed in \cite{BW1996}.
It is noticed in \cite{LiuGuoWu5} that $\D$ (hence, $\D^c$)  is an $f$-simplicial complex.
Then we take advantage of Eagon-Reiner theorem (\cite[Theorem 8.1.9]{Herzog and Hibi}) and CoCoA (\cite{CoCoA2020}) to check that both  simplicial complex $\D$  and its complement $\D^c$ are minimal CM. In particular, neither $\D$ nor $\D^c$ is shellable, which is hard to check without the notion of minimal CM (refer to \cite[Example 7.7]{BW1996} for a general theoretical treatment).
}
\end{myeg}

Note that {\it a permutation on the set $[6]$ may produce a new
simplicial complex $\D_1$}, which has the same property with $\D$. For example, the permutation $(1,2,3,4,5,6)$ acts on $\mathcal F(\D)$ and produces
$$\D_1=\langle \,234, 236, 124, 256, 125, 345, 135, 136, 456, 146\,\rangle.$$ Clearly, both $\D_1$  and $\D_1^c$ are $f$-simplicial complexes and minimal CM.

\section{3. Pure shellable versus Cohen-Macaulay}

We begin with the following  immediate consequence of Theorem \ref{homogeneous $f$-Ideal}:

\begin{mycor}\label{depth of f-complex} Let $\D$ be a $(d-1)$-dimensional pure $f$-simplicial complex with vertex set $[n]$. Then we have
$${\rm depth}(\D)=\begin{cases} d, \,\,\, \mbox{if } \D \,\, is \,\,CM \\ d-1, \hs{0.9cm} \mbox{if } \D \,\, is \,\,not\,\,CM. \end{cases} $$

\end{mycor}

\Proof Assume that $\D$ is not CM. Since $\D$ is a $(d-1)$-dimensional pure $f$-simplicial complex, $\mathcal F(\D)$ is an L-set, thus $\D^{(d-2)}=[n]_{d-1}$ holds. Clearly, $[n]_{d-1}$ is pure shellable, thus it is CM. Then the result follows from the fact that $${\rm depth}  (\D)= 1+{\max}\{ i\mid {\rm the \,\,i'th\,\,skeleton}\,\,\D^{(i)}\,\, {\rm is \,\,CM.}   \}.\quad\quad\qed$$

If the minimal subcomplex $\Gamma $ in Theorem \ref{shelled over main} is a simplex, say,  $\G=\langle F_0\rangle$, then it follows by Lemma \ref{shelling move} that the following is a shelling of $\D$, thus $\D$ is shellable:
$$F_0,F_1,F_2,\ldots,F_j.$$
To be more precisely, we have

\begin{mthm}\label{shellable via CM} For a pure simplicial complex $\D$, the following statements are equivalent:

$(1)$ $\D$ is shellable.

$(2)$ There exists a full sequence of subcomplexes $\D_i$ such that all $\D_i$ are CM, i.e., there is a total order $F_j,F_{j-1}, \ldots, F_1, F_0$ of all facets of $\D$ such that each $\D_i=:\langle F_0,F_1,\ldots,F_i\rangle $ is CM for $j\geq i\geq 1$, or equivalently,
each ideal $I(\D_i^c)$ has a linear resolution.\index{simplicial complex!pure shellable}
\index{simplicial complex!Cohen-Macaulay}\index{simplicial complex!subcomplex}

\end{mthm}

\Proof  $(1)\Longrightarrow (2):$  Let $F_0,F_1,\ldots, F_j$ be a shelling of $\D$ and let $\D_i=\langle F_0,F_1,\ldots, F_i\rangle$. Then for any $i$ with $1\le i\le j$, $\D_i=:\langle F_0,F_1,\ldots,F_i\rangle$ is pure and shellable, thus is CM.

$(2)\Longrightarrow (1):$ Let $F_j,F_{j-1}, \ldots, F_1, F_0$ be a full sequence of facets of $\D$ such that each $\D_i=:\langle F_0,F_1,\ldots,F_i\rangle $ is CM for $j\geq i\geq 1$.
  Then by Lemma \ref{shelling move},
  $$\D_{j-1}\cap \langle F_j\rangle, \D_{j-2}\cap \langle F_{j-1}\rangle,\ldots \D_{1}\cap \langle F_1\rangle$$
  are all pure of dimension $\dim \D-1$. By definition, $F_0,F_1,\ldots , F_j$ is a shelling of $\D$, thus $\D$ is a shellable simplicial complex. The rest statement follows from Eagon-Reiner theorem, and is convenient for checking by applying CoCoA.\quad\quad\qed

\vs{2mm}Clearly, Theorem \ref{shellable via CM} shows the exact relation between the conditions of shellable and CM for a pure simplicial complex. It also exhibits the importance of Lemma \ref{shelling move}.

As is well-known, it is in general a hard work to check  if a pure simplicial complex is  shellable. It seems that the new concept shelled over could open an algorithmic gate on attacking this problem, based on the algebraic characterization of a CM simplicial complex by Eagon-Reiner theorem. Refer to Examples \ref{f-idealNonShellable}, $3.4$ and
\ref{[8]4} for concrete operations and calculations.

When considering the condition of connected in codimension 1 (\cite[Proposition 1.12]{ZhengXX}), we have the following easy observation:

\begin{myprn} Let $\D$ be a pure simplicial complex of dimension $d-1$, which is not a simplex. Consider the following conditions:

$(1)$ For each face $\si$ of $\D$ such that ${\rm dim}\,{\rm lk}_\D(\si)>0$, ${\rm lk}_\D(\si)$ is connected.

$(2)$ For each facet $F$ of $\D$, $\D$ is a shelling move  of $\D_F$.

\noindent Then $(1)$ implies $(2)$.

\end{myprn}

\Proof  This follows from the proof of \cite[Lemma 3.1]{Dao2020}.\quad\quad\qed

\vs{2mm}The converse does not hold true in general. For example, the simplicial complex $\D=:\langle 1234,1235, 1278,1279\rangle$ is a shelling move over $\D_F$ for each facet $F$, but ${\rm lk}_\D(12)=\langle34,35,78,79\rangle$ and it is  disconnected.

For a simplicial complex $\D$, recall from \cite[Lemma 1.5.3]{Herzog and Hibi} that $I_{\D^\vee}=I(\D^c)$ holds, where $\D^c=\langle \, [n]\smallsetminus F\mid F\in \mathcal{F}(\D)\,\rangle$ and, $\D^\vee$ is the Alexander dual complex of $\D$.  Recall that $\D$ is said to be CM, if the Stanley-Reisner ideal $I_{\D}$ of $\D$ is CM. Recall also the Eagon-Reiner theorem (\cite[Theorem 8.1.9]{Herzog and Hibi}), i.e., a simplicial complex  $\D$ is CM
\iff the Stanley-Reisner ideal $I_{\D^\vee}$ of $\D^\vee$  has linear resolution. Thus $\D$ is CM \iff the monomial ideal $I(\D^c)$ has linear resolution. These results together with CoCoA are crucial to our next work.

In the following, we consider
simplicial complexes of   kind $(8,4)$
and apply Theorem \ref{shellable via CM} and Example \ref{f-idealNonShellable} to the following construction:

\begin{myeg}\label{shelledoverExg} {\rm  We start from the set $$A=\{1345, 1347,  1358, 1367, 1368, 1456, 1468, 1478, 1567, 1578\}.$$ By Example \ref{f-idealNonShellable}, it generates a minimal CM simplicial complex $\Gamma$, where $$\Gamma=\langle\, \{F\mid F\in A\}\,\rangle .$$

$(1)$ Let
$$B=\{1234, 1235, 1246, 1258, 1357, 1458, 1568, 2345, 2347, 2346, 2356, 2457,2468, $$$$
 2578, 2678, 2467,   2456, 2567, 2367,3456, 3478, 3678, 4567, 4578, 4678\}, $$
and let $\D_1=\langle \{F\mid F\in A\cup B\rangle.$
Then it is checked that $\D_1$ is a CM simplicial complex via CoCoA (\cite{CoCoA2020}), and the following sequence
$$F_{25}=1246, 1258, 1235, 1357, 3478, 1568, 1234, 2345, 2356, 2457, 2468, 2578, 2678,$$
$$ 2567, 2456, 3678, 3456, 2346, 1458, 4578, 4567, 2367, 2467, F_2=4678, F_1=2347$$
of facets are found to make the CM simplicial complex $\D_1$ shelled over the minimal CM simplicial complex $\Gamma$, where each of the $24$ monomial ideals
$$I(\langle \mathcal(\D_1^c)\smallsetminus \{ F_{25}^c\}\rangle),
I(\langle \mathcal F(\D_1^c)\smallsetminus \{F_{25}^c, F_{24}^c\}\rangle),\ldots, I(\langle \mathcal F(\D_1^c)\smallsetminus \{F_{25}^c, \ldots,F_{2}^c\}\rangle)$$
is tested via CoCoA to have linear resolution.

Note that $\D_1$
is not an $f$-simplicial complex since $127$ is not in the lower set of $\mathcal F(\D_1)$, i.e., $\mathcal F(\D_1)$ is not an $L$-set. We do not know if $\D_1$  is shellable.

$(2)$ Inspired by the previous construction, we now construct an  $f$-simplicial complex
$$\D_2=\langle\, \{F\mid F\in A\cup C\}\,\rangle,$$ where $$C=\{1235,   1236, 1237, 1247, 1268,     1258,    1357, 1458,   1568,  2345, 2347,  2348,  2356, $$
 $$ 2457, 2467, 2468, 2578, 2678,
  3456, 3467,  3478, 3678, 4567, 4578, 4678\} .$$
We checked that $A\cup C$ is an LU-set, thus $\D_2$ is indeed an $f$-simplicial complex. We checked that $\D_2$ is CM via CoCoA. Furthermore, we claim that the complex $\D_2$ is shellable with the shelling $F_0,F_1,\ldots,F_{34}$, where
$$F_{34}=1247,F_{33}=1237,1357,1358,1368, 2347, 2348, 2356,1236, 1235,$$
$$1268,1258, 1345, 1347,1367, 1567, 1568, 1578, 1458, 1456, 1468, 1478, 3478,$$
$$ 2345,2457, 2467,  2468, 2578, 2678,  3456,3467, 3678,  4567,   F_1=4678,   F_0=4578.$$
In fact, we use CoCoA to show that each $I(\langle F_0^c,\ldots, F_i^c\rangle)$ has linear resolution ($\forall 34\ge i\ge 1$), thus all
$$\langle F_0, F_1,\ldots,F_i \rangle $$
are Cohen-Macaulay simplicial complexes. Then it follows from  Theorem \ref{shellable via CM} that $\D_2$ is shellable with
$$F_0,F_1,\ldots,F_{34}$$
as a shelling.

It is natural to ask if $\D_2$ is shelled over $\Gamma$ constructed in $(1)$? We tried this via CoCoA, and the answer is yes. The following sequence
$$G_{25}=1247, G_{24}=1357, 1237, 2347, 2348, 2345, 2356,  1235, 1236, 1268, 1258,2578,$$
$$ 2678, 2468, 2467, 2457, 3478,3678,1568,4678,3456,3467,4567,G_2=4578, G_1=1458$$
of facets are found to make the CM simplicial complex $\D_2$ shelled over the minimal CM simplicial complex $\Gamma$, where all the $24$ facet monomial ideals $I(\langle \mathcal F(\D_2^c)\smallsetminus \{G_{25}^c, \ldots,G_{r}^c\}\rangle)$ are tested to have linear free resolutions.

Note that for the same CM simplicial complex $\D_2$, the first minimal CM subcomplex is $\langle F_0 \rangle$ and it has only one facet, while the second minimal CM subcomplex is $\G$ and it has ten facets. This is the end of Example \ref{shelledoverExg}.
}
\end{myeg}

Finally, note that Example \ref{f-idealNonShellable} provides a {\em very well-distributed} simplicial complex, i.e., each number $r$ in $[6]$ appears 5 times in the facets. Motivated by this observation, we now construct a very well-distributed simplicial complex whose facet set contains 34 elements in $[8]_4$.

\begin{myeg} \label{[8]4}{\rm Let $\mathcal F(\D^c)$ be the set $D$ consisting of $34$ elements, where
$$D=\{ 1234, 1235,
1246, 1247,
1258, 1345,
1358, 1367,
1368, 1378,
1456, $$
$$
1457,1467, 1478,
1568, 1578,
1678, 2346,
2348, 2356,
2358, 2367,
2378, $$
$$2456,2467, 2468,
2478, 2567,
2678, 3456,
3457,  3568,
3578,  4578\}.$$
We checked the following:

$(1)$ $D$ is very well-distributed, i.e., it has type  $1^{17}2^{17}3^{17}4^{17}5^{17}6^{17}7^{17}8^{17}$.

$(2)$ $D$ is an LU-set over $[8]$, so that adding any element from $[8]_4\smallsetminus D$ can generate a pure  $f$-simplicial complex $\Gamma^c$.

$(3)$ The ideal $I(\D^c)$ has the following linear resolution, thus $\D$ is not CM:
{\fs $$0\LRA R(-8)^2\LRA R(-7)^{21}\oplus R(-8)\LRA R(-6)^{68}\oplus R[-7]^2\LRA R(-5)^{81}\oplus R[-6]\LRA R(-4)^{34}\LRA R.$$}

Note that $\D^c$ has the same properties.

Among the $36$ $f$-simplicial complexes $\Gamma$ obtained in $(2)$, $7$ are CM. In fact,  $F_{35}^c$ can be chosen as anyone of the following:
$$1236,1238,1268,1346,1348,1468,2368$$
such that $I(\langle D^c\cup \{F_{35}^c\}\rangle)$ has linear resolution.
Note that unfortunately, none of the $7$ CM simplicial complexes is minimal CM. Furthermore, all $\Gamma$ are shelling moves of $\Gamma_{F_2}$, where $F_2^c=1235$.

Finally, we consider $D'=:[8]_4\smallsetminus D$, which consists of $36$ elements, as follows:
$$D'=\{1236,1237,1238,1245, 1248, 1256,1257,1267, 1268,1278, 1346, 1347, \}$$
$$1348,1356, 1357, 1458, 1468,1567, 2345, 2347, 2357, 2368, 2457, 2458, $$
$$2568, 2578,3458, 3467, 3468, 3478, 3567, 3678, 4567, 4568, 4678,5678\}.$$
Deleting any element will result in a homogeneous complement of some $\Gamma^c$ in $(2)$. We use CoCoA to calculate the $36$ $I((\Gamma^c)')$ and, find $15$ CM simplicial complexes $\Gamma'$. The following are all elements when one of which is deleted, the corresponding $I((\Gamma^c)')$ has linear resolution, thus $\Gamma'$ is CM:
$$1238, 1257, 1267, 1356, 1357, 1458,1567, 2345, 2357, 2578,3567, 3678, 4567, 4678, 5678.$$
Note that $(\Gamma^c)'=(\Gamma')^c$ always holds true. It also gives a negative  answer to Question \ref{Homogeneous complement question}.
This is the end of the example.
}
\end{myeg}

Note that in many examples of CM simplicial complexes $\D$, we have $\D_F\cap \langle F\rangle =\langle \,G\mid G\in\partial F\,\rangle $ holds true
for most of the facets $F$, but not in all cases, as the following example shows:

\begin{myeg} \label{[8]41} Let $\mathcal F(\D_3^c)=D\cup \{1236\}$, in which $D$ is taken as in Example \ref{[8]4}. Let $F^c=4578$ and consider $\Gamma=:\D_F\cap \langle F\rangle $. Then we have $123\not\in \mathcal F(\Gamma) $, thus $\mathcal F(\Gamma)$ is a proper subset of $\partial F=:\{123,126,136,236\}$.
\end{myeg}

Notice the following
$$30=5\times 6=10\times 3,\,\, 34\times 4=8\times 17=136<140=35\times 4.$$
Notice the following fact:

$(1)$ For any odd number $d\ge 3$, $d\times \frac{1}{2}\binom{2d}{d}$  is divided by $2d$, i.e., $4\mid \binom{2d}{d}$ holds true.

\noindent Based on the examples and Theorem \ref{homogeneous $f$-Ideal}, we now pose the following:

\begin{myconj} $(a)$ There exist in $[2d]_d$ very well-distributed $f$-simplicial complexes which are minimal CM, if one of the following conditions holds true:

\hspace{1cm} $(1)$  $d\ge 3$ and $d$ is an odd number.

\hspace{1cm}$(2)$ $d\ge 4$ and $d$ is an even number  such that $4\mid \binom{2d}{d}$.

$(b)$ In $[8]_4$, there  exists no  very well-distributed $f$-simplicial complex which is minimal CM.
\end{myconj}

Note that $4\nmid\binom{8}{16}$ also holds. Thus in $[16]_8$, the pure simplicial complexes generated by $\frac{1}{2}\binom{8}{16}$ subsets may perhaps behave just like the pure simplicial complexes in $[8]_4$.

\section{4.  Minimal Cohen-Macaulay $f$-simplicial complexes}

In this section, we study properties of  $f$-simplicial complexes which are minimal CM.  For this, we need the following:

\begin{mylem}\label{kernal of chain complex of a simplex} Let $$\mathfrak C: 0\LRA C_{n-1}\overset{\partial_{n-1}}{\LRA} C_{n-2}\overset{\partial_{n-2}}{\LRA} \cdots\overset{\partial_2}{\LRA} C_1\overset{\partial_1}{\LRA} C_{0}\LRA 0$$ be the chain complex of the simplex $\langle\,[n]\,\rangle$ over a field $\mathfrak{K}$. Then the $\mathfrak{K}$-subspace $\ker(\partial_r)$ has dimension $\binom{n-1}{r+1}$.
\end{mylem}

\Proof For $1<r\leq n-1$, let $\sum_{1\leq i_1<i_2<\cdots<i_r\leq n}x_{i_1i_2\ldots i_r}\cdot i_1i_2\ldots i_r\in {\rm ker}\,\partial_{r-1}$, where the second $i_1i_2\ldots i_r$ denotes the subset $\{i_1,i_2,\ldots, i_r\}$ of $[n]$ and $x_{i_1i_2\ldots i_r}$ are elements of the base field $\mathfrak K$. Then
$$\sum_{1\leq i_1<i_2<\cdots<i_r\leq n}x_{i_1i_2\ldots i_r}\sum_{j=1}^r (-1)^{j-1}i_1i_2\ldots i_{j-1}\hat{i_j}i_{j+1}\ldots i_r=0$$
holds true. Since $C_{r-2}$ is a free $\mathfrak K$-module with basis $[n]_{r-1}$, we have got a system of homogeneous linear equations, which consists of $\binom{n}{r-1}$ equations with $\binom{n}{r}$ variable $x_{i_1i_2\ldots i_r}$. We write these $x_{i_1i_2\ldots i_r}$ as well as $i_1i_2\ldots i_{r-1}$ in lexicographic order, and consider the rank of the coefficient matrix $M_{\binom{n}{r-1}\times \binom{n}{r}}$. Clearly, the first $\binom{n-1}{r-1}$ columns, i.e., the coefficients of $x_{1i_2\ldots i_r}$, are linearly independent. It can be checked that each other column is a linear combination of them. Furthermore, for $i_1i_2\ldots i_r$ with $1\not\in \{i_1,\ldots,i_r\}$,
note that $$\partial_{r-1}(i_1i_2\ldots i_r)=i_2\ldots i_r-i_1i_3\ldots i_r+\cdots+(-1)^{r-1}i_1\ldots i_{r-1},$$
in the $i_1\ldots i_r$-th column vector $v_{i_1\ldots i_r}$ of $M$, the $i_1\ldots \hat{i_{j}}\ldots i_r$-th component is $(-1)^{j-1}$ ($1\leq j\leq r$) and, all other components are zero. Thus we have
\begin{equation}\label{relation}v_{i_1\ldots i_r}=v_{1i_2\ldots i_r}-v_{1i_1i_3\ldots i_r}+\cdots +(-1)^{r-1}v_{1i_1\ldots i_{r-1}}.
\end{equation}
The exact details are essentially the same with the verification of the fact that a cone is acyclic, refer to
\cite[Proposition 5.2.5]{Villarreal}.

Finally,  we proved that the dimension of the vector space ${\rm ker}\,\partial_{r-1}$ is the following
$$\binom{n}{r}-\binom{n-1}{r-1}=\binom{n-1}{r}.$$
Note that the key to calculate the kernel of general $\partial_{i-1}$ is the equality $$\binom{n}{i}=\binom{n-1}{i}+\binom{n-1}{i-1}.$$ This is the end of the verification.
\qed

{\bf Remark.} We illustrate the proof in computational way in two particular  cases.

The first case is $[6]_3$, and we check that
$\ker\,\partial_1$ has dimension $\binom{n-1}{2}$, where $n=6$. In fact, let $\sum_{1\leq i<j\leq 6}x_{ij}\{i,j\}\in {\rm ker}\,\partial_1$, we have $$0=\sum_{1\leq i<j\leq 6}x_{ij}(\{j\}-\{ i\}).$$ Since $C_0$ in the chain complex is a free $\mathfrak{K}$-module with basis $\{1\},\{2\},\ldots,\{6\}$, we get the following system of linear equations:
$$\begin{cases}
\sum_{i=1}^6-x_{1i}=0 \\
x_{12}-\sum_{i=2}^6x_{2i}=0\\
x_{13}+x_{23}-\sum_{i=3}^6x_{3i}=0\\
\sum_{i=1}^3x_{i4}-x_{45}-x_{46}=0\\
\sum_{i=1}^4x_{i5}-x_{56}=0\\
\sum_{i=1}^5x_{i6}=0.
\end{cases} $$

The coefficient matrix is
 \begin{displaymath}
 \left( \begin{array}{cccccccccccccccc}
-1&-1&-1&-1&-1&&&&&&&&&&\\
1 &&&&&-1 &-1 &-1 &-1 &&&&&&\\
 &1 && &&1&&&&-1 &-1 &-1 &&& \\
 & &1 & & &&1&  &&1&&&-1 &-1 & \\
&&&1 &   &&&1&&  &1&&1&&-1 \\
 & & & &1 & &  & &1&& &1& &1&1
\end{array} \right),
\end{displaymath}
\noindent hence, the solution of the system has exactly
$$10=\binom{5}{2}=\binom{6}{2}-\binom{5}{1}$$
free variables, they are all $k_{ij}$ except these $\{i,j\}$ including 1.
Note that for a general $n$,  $\dim\, \ker(\partial_1)=\binom{n-1}{2}$ is verified in an exactly same way. Note also that in the coefficient matrix, we have column vector relation
$v_{23}=v_{12}-v_{13}.$

\begin{center}
\includegraphics[width=15.8cm]{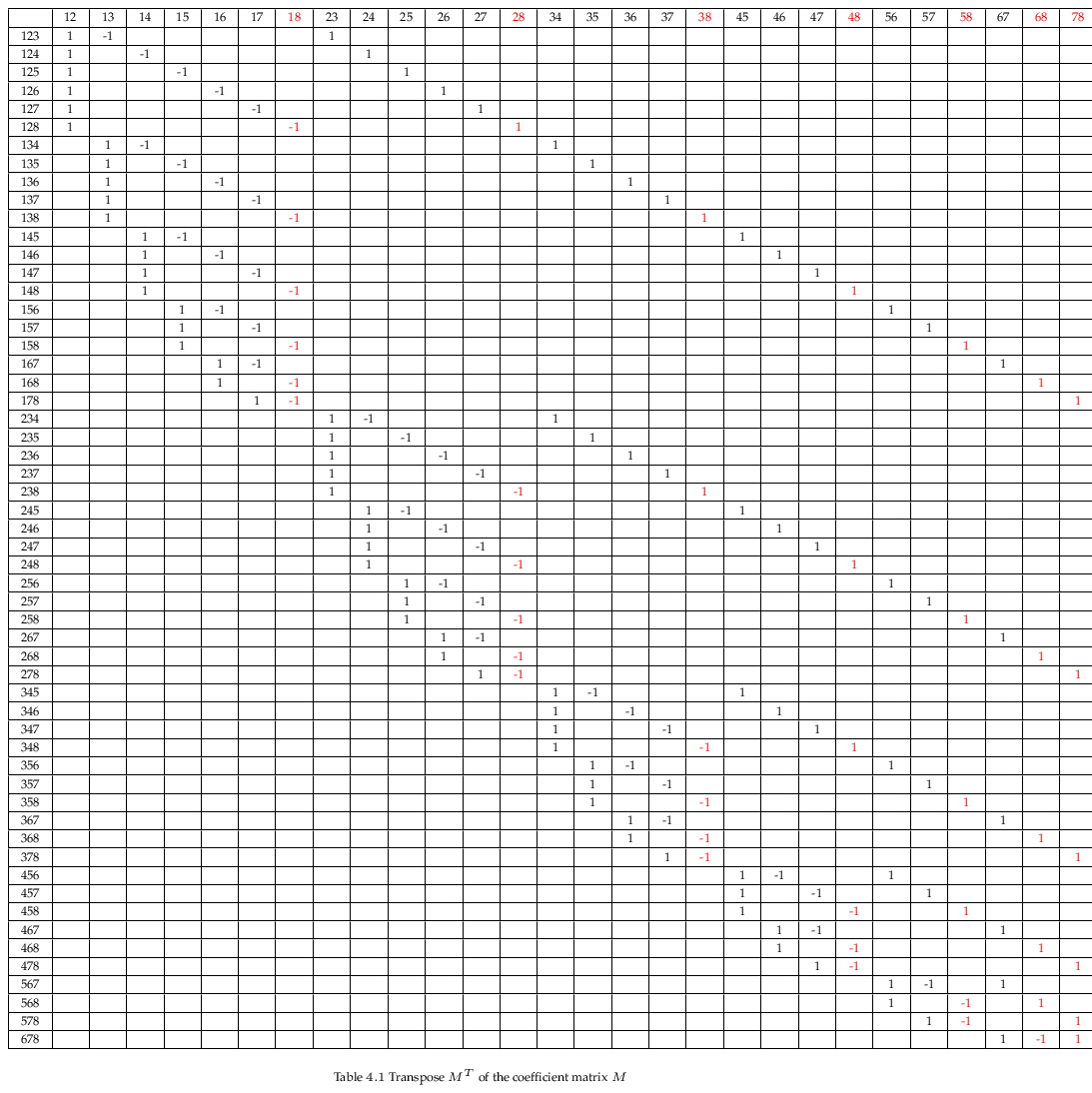}
\end{center}

The second case is $[8]_4$, and we check that
$\ker\,\partial_2$ has dimension $\binom{n-1}{3}$, where $n=8$. In fact, let $\sum_{1\leq i<j<k\leq 8}x_{ijk}\{i,j,k\}\in {\rm ker}\,\partial_2$, we have $$0=\sum_{1\leq i<j<k\leq 8}x_{ijk}(\{j,k\}-\{i, k\}+\{i,j\}).$$
Since $C_1$ in the chain complex is a free $\mathfrak{K}$-module with basis set $\{ij\mid 1\leq i<j\leq 8\},$ we get a system of linear
equations, which consists of $\binom{8}{2}$ linear equations with $\binom{8}{3}$ variables $x_{ijk}$. We write out the coefficient matrix $M$ in Table $4.1$ and it is clear that the row (in the table) rank of the matrix is not less than $\binom{7}{2}$. Actually,  after doing Gaussian elimination via excel, it is calculated that the matrix
 has rank $\binom{7}{2}=:21$. We also checked the rank by taking advantage of CoCoA (\cite{CoCoA2020}) $5.3.3.$ Besides,
all $x_{ijk}$ except $x_{1jk}$'s can be chosen as free variables. Certainly, there is an alternative explanation as appeared in Table $4.1$. This shows $$\dim\,\ker(\partial_2)=\binom{8}{3}-\binom{7}{2}=\binom{7}{3},$$
as is claimed. Note that for a general $n$,  $\dim\, \ker(\partial_2)=\binom{n-1}{3}$ is verified in a completely same way.

Note also that in Table $4.1$, we have row vector relation
$v_{234}^T=v_{123}^T-v_{124}^T+v_{134}^T,$
which is a particular case of (\ref{relation}).

\begin{mycor} \label{homologies of $[n]_r$} Let $\D=\langle\, \{\si\mid \si\in [n]_r\}\,\rangle $. Then we have

$${\rm dim}\,\tilde{H_i}(\D)=\begin{cases}  \binom{n-1}{r}, \,\,\, \mbox{if } \,\,i=n-1 \\ 0, \hs{0.9cm} \mbox{otherwise }. \end{cases} $$


\end{mycor}

\begin{mthm} \label{acyclic3} Let $\D$ be a simplicial complex over vertex set $[n]$ with $\dim\,\D=d-1\ge 0$. If $\D$ is an $f$-simplicial complex and it is minimal CM, then $d\ge 3$ and $n=2d$.

\index{simplicial complex!minimal CM}\index{simplicial complex!acyclic}
\end{mthm}

\Proof   It is known that connected simplicial complexes of dimension $1$ are shellable and pure. On the other hand, if $\D$ is not connected, then it is not CM. So, there exists no minimal CM $f$-simplicial complexes of dimension 1.

Now let $\D$ be an $f$-simplicial complex of dimension $d-1$, which is minimal  CM.
Then
$\tilde{H}_{d-1}(\D)=0$ by Theorem \ref{facetdepth},  which means that $\partial_{d-1} $ is injective.  Since $\D$ is an $f$-simplicial complex, $\mathcal F(\D)$ is an L-set, hence
$\D^{(d-2)}=[n]_{d-1}$ holds true, thus, $\tilde{H}_{i}(\D)=0, \forall 0\leq i\leq d-3$. Hence $\D$ is acyclic \iff $\tilde{H}_{d-2}(\D)=0$, and the latter holds true \iff ${\rm dim}_{\mathfrak{K}}\,{\rm ker}(\partial_{d-2})=\frac{1}{2}\binom{n}{d}$ by Theorems \ref{facetdepth} and \ref{homogeneous $f$-Ideal}.

By Lemma \ref{kernal of chain complex of a simplex}, we have ${\rm dim}_{\mathfrak{K}}\,{\rm ker}(\partial_{d-2})= \binom{n-1}{d-1} $.
We get
$$\frac{n(n-1)\cdots(n-d+1)}{2\cdot d!}= \frac{(n-1)(n-2)\cdots(n-d+1)}{(d-1)!},$$
since a minimal CM simplicial complex is always acyclic by \cite{Dao2020}. Thus we have $n=2d$.
\quad\quad\qed

By Mayer-Vietoris long exact sequence theorem, we get\index{Mayer-Vietoris Long Exact Sequence}

\begin{mycor} If $\D$ is an $f$-simplicial complex generated by a subset of $[2d]_d$ and it is minimal CM, then $\tilde{ H}_i(\D_F)\cong \tilde{ H}_i(\langle F\rangle\cap \D_F)$  holds for all integer $i$, where $F$ is any facet of $\D$ and, $\langle F\rangle\cap \D_F$ is pure of dimension $d-2$.

\end{mycor}

\Proof For any facet $F$ of $\D$, let
$$\D_1=\D_F, \, \D_2=\langle\, F\,\rangle,\, \D_3=\D_1\cap \D_2 =:\D_F\cap \langle\, F\,\rangle.$$
Then $\D=\D_1\cup \D_2$. By \cite[Theorem 25.1,  page 142]{Munkres}, we have the following long exact sequence of $\mathfrak{K}$-spaces:
$$0\longrightarrow \tilde{H}_{d-1}(\D_3)\longrightarrow  \tilde{H}_{d-1}(\D_1)\oplus \tilde{H}_{d-1}(\D_2)\longrightarrow  \tilde{H}_{d-1}(\D)\hspace{2cm}$$
$$\hspace{-1.8cm} \overset{\partial_{d-1}}{\longrightarrow} \tilde{H}_{d-2}(\D_3)\longrightarrow \tilde{H}_{d-2}(\D_1)\oplus \tilde{H}_{d-2}(\D_2)\longrightarrow  \tilde{H}_{d-2}(\D)$$
$$\hspace{-9cm}\overset{\partial_{d-2}}{\longrightarrow}\,\,\,  \cdots\cdots\cdots$$
$$\hspace{-3.5cm}\overset{\partial_2}{\longrightarrow} H_{1}(\D_3)\longrightarrow  H_{1}(\D_1)\oplus H_{1}(\D_2)\longrightarrow  H_{1}(\D)$$
$$\hspace{-2.3cm}\overset{\partial_1}{\longrightarrow} \tilde{H}_{0}(\D_3)\longrightarrow  \tilde{H}_{0}(\D_1)\oplus \tilde{H}_{0}(\D_2)\longrightarrow  \tilde{H}_{0}(\D)\longrightarrow0.$$
Note that Lemma \ref{shelling move} implies that $\D_3$ is pure of dimension $d-2$, while it follows from Theorem \ref{acyclic3} that $\tilde{H}_i(\D)=0$ holds for all $i$, and $\tilde{H}_i(\D_2)=0$ holds true clearly. Then
$$0\longrightarrow \tilde{H}_i(\D_3)\longrightarrow \tilde{H}_i(\D_1)\longrightarrow 0$$
is an exact sequence for every $i$.\quad\quad\qed

\end{document}